\documentclass[12pt]{article}
\usepackage{fullpage}

\begin{document}

\title{Evaluation of Authors and Journals}

\author{Joseph B. Keller\\Departments of Mathematics and Mechanical Engineering\\Stanford University\\Stanford, CA  94305}

\date{October 1985}

\maketitle
\begin{abstract}
A method is presented for evaluating authors on the basis of citations.  
It assigns to each author a citation score which depends upon the number of times he is cited, and upon the scores of the citers. 
The scores are found to be the components of an eigenvector of a normalized citation matrix.
The same method can be applied to citation of journals by other journals, to evaluating teams in a league \cite{keller}, etc.

\end{abstract}

\section{Introduction}

One commonly used measure of the influence of an author is the number of times his work is referred to by others in a given period of time.
For a scientific author, this number can be found by counting the number of citations of his work listed in the Science Citation Index for that period.
However, this measure fails to take into account who is doing the citing.  A citation by an influential author ought to carry more weight than one by an 
unknown author, but a simple count of citations does not give it more weight.  A more appropriate measure would be a weighted count of the citations,
in which each citation is weighted by some measure of the influence of the citer.  We shall show how to find such a measure, which we call a \emph{citation score\/}, 
or just \emph{score} for short.

We begin by assuming that each author can be assigned a score, and that the weight of a citation is proportional to the score of the citer.
Then we determine the score of an author by adding up the weights of all the citations of his work.  The circularity of this method leads to a requirement
of consistency among the scores, which determines them as the solutions of a system of linear equations.  This method also applies to the evaluation
of journals on the basis of citations of them in other journals.  

We have used this methods before \cite{keller} to evaluate teams in a league, with $c_{ij}$ the number of times
that team $i$ beats team $j$.  Other ranking methods are reviewed by Moon and Pullman \cite{mp}.

\section{The citation matrix}

Let us consider $N$ authors, numbered from $1$ to $N$.  We denote by $c^\prime_{ij}$ the number of times that $j$ cited $i$, omitting self citations,
so that $c^\prime_{ii}=0$.  The  $c^\prime_{ij}$ form an $N$ by $N$ square matrix $C^\prime$, which we call the citation matrix.

The $j$-th column of $C^\prime$ records all the citations by $j$ of others, and the sum of the entries in that column is the total number of citations by $j$.  If the sum is not zero, we
normalize the column by dividing each entry by the column sum.  Thus we define $c_{ij}$ by
\begin{equation}
c_{ij} =c^\prime_{ij} \biggr/  \sum\limits^N_{k=1} c^\prime_{kj}.
\label{1}
\end{equation}
If all the entries in column $j$ are zero, we define $c_{ij}=0$.  We call the matrix with entries $c_{ij}$ the normalized citation matrix $C$.
From the definition (\ref{1}) we see that $c_{ij}$ is the fraction of $j$'s citations which refer  to $i$.  Furthermore, each column sum is unity, unless all the entries in the column are 
zero, in which case it is zero.

Next we denote by $x_i$ the score of author $i$.  According to the method mentioned in the Introduction, $x_i$ is given by\
\begin{equation}
x_i = \lambda^{-1} \sum\limits^N_{j=1}\,  c_{ij} x_j.
\label{2}
\end{equation}
Here $\lambda^{-1}$ is a factor of proportionality.  Thus $x_i$ is the sum of contributions from all citers $j\not= i$.  Each contribution is the product of the 
score of $j$ times the fraction of $j$'s citations which refer to $i$, all times $\lambda^{-1}$.

Equation (\ref{2}) is a system of $N$ linear homogenous equations for the scores $x_i$.  In terms of the score
 vector $x=\left(x_1,\ldots,x_N\right)$, (\ref{2}) can be written
\begin{equation}
Cx=\lambda x.
\label{3}
\end{equation}
Thus the score vector $x$ is an eigenvector of the normalized citation matrix $C$ corresponding to the eigenvalue $\lambda$.

\section{Eigenvectors of the citation matrix}

In order for $x$ to be a score vector, its components must be non-negative numbers.  Since  all the entries $c_{ij}$ of $C$ are non-negative, Theorem 3 on p.\ 66 of Gantmacher \cite{gant} shows that $C$ has a real non-negative eigenvalue $\lambda$ with a non-negative eigenvector.

This non-negative eigenvector $x$ could be used to determine the scores if it were unique, aside from a constant factor.  It will certainly be unique if $C$ is irreducible, i.e., if it cannot be put in the following form by a permutation of the indices:
\begin{equation}
C=\left(\begin{array}{cc}
A & 0\\
B&D
\end{array}\right).
\label{4}
\end{equation}
Here $A$ and $D$ are square matrices.

When $C$ is irreducible, Frobenius' theorem states that it has a real non-negative eigenvalue $\lambda$ which is larger than the modulus of any other eigenvalue.  Furthermore the eigenvector $x$ corresponding to $\lambda$ is unique up to a scalar factor, and all its components are positive.  (Gantmacher, p.\ 53, theorem 2.)  In this case the components of $x$, with some suitable normalization, can be used as the scores.

The eigenvalue $\lambda$ can be determined by first summing (\ref{2}) over $i$.  Then from the fact that $C$ is normalized, it follows that the sum of $c_{ij}$ over $i$ is unity.  The other possibility, that all the $c_{ij}$ in one column vanish, cannot occur when $C$ is irreducible.  Thus we obtain from (\ref{2})
\begin{equation}
\sum\limits^N_{i=1} x_i =\lambda^{-1} \sum\limits^N_{j=1} \left(\sum\limits^N_{i=1} c_{ij}\right) x_j =\lambda^{-1} \sum\limits^N_{j=1} x_j.
\label{5}
\end{equation}
Now (\ref{5}) shows that $\lambda=1$.

Since $\lambda =1$, we can write (\ref{3}) in the form
\begin{equation}
Cx=x.
\label{6}
\end{equation} 
We also impose the normalization condition that the largest component of $x$ is unity;
\begin{equation}
\max\limits_{i} \, x_i =1.
\label{7}
\end{equation}
The previous considerations show that when $C$ is irreducible, (\ref{6}) and (\ref{7}) have a unique solution in which all the components $x_i$ are positive.  They are the scores we wanted.

\section{The reducible case}

Let us now suppose that the normalized citation matrix $C$ is reducible, and that it has been put in the form (\ref{4}). If $B=0$, $C$ is block diagonal, which means that there are two distinct sets of authors with no references by members of either set to the works of those in the other set.  This might be the case if the two sets of authors write about completely separate fields.  When $C$ is block diagonal, it has at least two linearly independent eigenvectors, one corresponding to $A$ and another to $D$.  It is not surprising that in this case the scores of the two sets of authors are unrelated.

Next we suppose that $B\not= 0$, and that both $A$ and $D$ are irreducible.  Then (\ref{5}) still holds, so $\lambda =1$.  Now we write $x$ in the partitioned form $x= (y,z)$, and (\ref{3}) becomes the pair of equations
\begin{eqnarray}
Ay=y, \label{8}\\
By+Dz=z.
\label{9} 
\end{eqnarray}

The matrix D-I in (\ref{9}) is singular because $\lambda=1$ is an eigenvalue of $D$.  Furthermore, $\omega$,
the corresponding left eigenvector of $D$, has all positive components.  The inner product of $\omega$ with (\ref{9}) yields the solvability condition $\omega^TBy=0$.  Since all the entries in $B$ are non-negative, while $\omega$ and $y$ are positive, $\omega^T By>0$ unless $y=0$.  Thus $y=0$ and then (\ref{9}) has a positive solution.  Therefore when $C$ is reducible to the form (\ref{4}) with $A$ and $D$ irreducible and $B\not= 0$, there is a unique normalized score vector $x= (y,z)$ in which $y=0$.

The conclusion that $y=0$ for all the authors corresponding to $A$ is unsatisfactory, and suggests that the reducible case should be treated differently.

\end{document}